\documentclass[11pt, reqno]{amsart}
\usepackage{amsfonts,latexsym}
\usepackage{amsmath}
\usepackage{amscd}
\usepackage{float,amsmath,amssymb,mathrsfs,bm,multirow,graphics}
\usepackage[dvips]{graphicx}
\usepackage[percent]{overpic}
\usepackage[numbers,sort&compress]{natbib}
\usepackage{enumerate}

\newcommand{\nequation}{\setcounter{equation}{0}}

\newcommand{\R}{{\Bbb R}}

\newcommand{\C}{{\Bbb C}}

\newcommand{\Z}{{\Bbb Z}}

\DeclareMathOperator{\im}{Im}
\DeclareMathOperator{\re}{Re}

\newcommand{\res}{\text{\upshape Res\,}}

\def\XXint#1#2#3{{\setbox0=\hbox{$#1{#2#3}{\int}$}
\vcenter{\hbox{$#2#3$}}\kern-.5\wd0}}

\allowdisplaybreaks

\newtheorem{figuretext}{Figure}

\input epsf
\title[The Unified Transform for the Modified Helmholtz]{\sc The Unified Transform for the Modified Helmholtz Equation in the Exterior of a Square}
\author{A. S. Fokas}
\address{A.S.F.: Department of Applied Mathematics and Theoretical Physics, University of Cambridge, Cambridge CB3 0WA, United Kingdom, and Research Center of Mathematics, Academy of Athens, 11527, Greece.}
\email{T.Fokas@damtp.cam.ac.uk} 

\author{J. Lenells}
\address{J.L.: Department of Mathematics and Center for Astrophysics, Space Physics \& Engineering Research, Baylor University, One Bear Place \#97328, Waco, TX 76798, USA.}
\email{Jonatan\_Lenells@baylor.edu}

\begin{document}

\begin{abstract} 
\noindent
The Unified Transform provides a novel method for analyzing boundary value problems for linear and for integrable nonlinear PDEs. The numerical implementation of this method to linear elliptic PDEs formulated in the {\it interior} of a polygon has been investigated by several authors (see the article by Iserles, Smitheman, and one of the authors in this book). Here, we show that the Unified Transform also yields a novel numerical technique for computing the solution of linear elliptic PDEs in the {\it exterior} of a polygon. One of the advantages of this new technique is that it actually yields directly the scattering amplitude. Details are presented for the modified Helmholtz equation in the exterior of a square.
\end{abstract}

\maketitle



%\tableofcontents

\section{Introduction}\nequation
The so-called Unified Transform, introduced by the first author in the late nineties \cite{F1997, Fbook}, has been implemented to a variety of boundary value problems for linear elliptic PDEs. In particular, it has led to the emergence of a new numerical technique for the Laplace, the Helmholtz, and the modified Helmholtz equations formulated in the {\it interior} of a polygon \cite{A2012, A2013, F2001, FFSS2009, FAS2013, FF2011, FFX2004, SSP2012, SFFS2008, SFPS2009, SPS2007, SSF2010}. This technique is illustrated for the particular case of a square in this book \cite{FISbookpaper}. A novel numerical implementation of the Unified Transform to the Helmholtz and modified Helmholtz equations formulated in the {\it exterior} of a polygon is presented in \cite{FLhelmholtz}. Here, we illustrate this implementation in the particular case of a square. 

\section{The Unified Transform}\nequation
The starting point of the Unified Transform is rewriting the given PDE as the condition that a one-parameter family of differential forms is closed. For the modified Helmholtz equation
\begin{align}\label{2.1}
  \frac{\partial^2 u}{\partial x^2} + \frac{\partial^2 u}{\partial y^2} - 4\beta^2 u = 0, \qquad \beta > 0, \quad (x,y) \in D \subset \R^2,
\end{align}
the relevant differential form is given by
\begin{align}\label{2.2}
  W(z,\bar{z}, \lambda) = e^{-i\beta(\lambda z - \frac{\bar{z}}{\lambda})} \left[\biggl(\frac{\partial u}{\partial z} + i\beta \lambda u\biggr) dz - \biggl(\frac{\partial u}{\partial \bar{z}} + \frac{\beta u}{i \lambda}\biggr)d\bar{z}\right],
\end{align}
where
\begin{align}\label{2.3}
z = x + iy, \qquad \bar{z} = x - iy, \qquad z \in D.
\end{align}
Indeed, recall that the differential of the one-form $W$,
$$W(z, \bar{z}) = A(z, \bar{z}) dz + B(z, \bar{z}) d\bar{z}, $$
is given by
\begin{align}\label{2.4}
  dW = \bigg(\frac{\partial B}{\partial z} - \frac{\partial A}{\partial \bar{z}}\bigg) dz \wedge d\bar{z}.
\end{align}
Using the above expression for the differential form $W$ defined by (\ref{2.2}) we find
\begin{align}\label{2.5}
dW = -2e^{-i\beta(\lambda z - \frac{\bar{z}}{\lambda})} \bigg(\frac{\partial^2 u}{\partial z \partial \bar{z}} - \beta^2 u\bigg) dz \wedge d\bar{z}, \qquad z \in D.
\end{align}
However, employing the change of variables defined by (\ref{2.3}), we find that the modified Helmholtz equation (\ref{2.1}) can be written in the form
\begin{align}\label{2.6}
\frac{\partial^2 u}{\partial z \partial \bar{z}} - \beta^2 u = 0, \qquad z \in D,
\end{align}
thus (\ref{2.5}) implies
\begin{align}\label{2.7}
dW = 0, \qquad z \in D.
\end{align}

\section{The Global Relation for Exterior Domains}\nequation
Suppose that the domain $D$ is exterior to some bounded domain $D_B$. Then equation (\ref{2.7}) and the complex form of Green's theorem imply
\begin{align}\label{3.1}
\int_{\partial D_B} W(z, \bar{z}, \lambda) = I(\lambda), \qquad \lambda \in \C,
\end{align}
where $I(\lambda)$ is defined by 
\begin{align}\label{3.2}
I(\lambda) = \lim_{r \to \infty} \int_{|z| = r} W(z, \bar{z}, \lambda), \qquad \lambda \in \C.
\end{align}
Assume that $u$ satisfies the Sommerfeld radiation condition, i.e.
\begin{align}\label{3.3}
\lim_{r \to \infty} \sqrt{r} \bigg(\frac{\partial}{\partial r} + 2\beta\bigg) u = 0.
\end{align}
It is well known (see \cite{K1961}) that as $r\to \infty$, $u$ admits the asymptotic expansion 
\begin{align}\label{3.4}
u \sim \sqrt{\frac{1}{\pi i \beta r}} e^{i(2i\beta r - \frac{\pi}{4})} \bigg(f_0(\varphi) + \frac{f_1(\varphi)}{r} + \frac{f_2(\varphi)}{r^2} + \cdots\bigg), \qquad r \to \infty.
\end{align}

Employing the above representation in (\ref{3.2}), it is shown in \cite{FLhelmholtz} that $I(\lambda)$ is given by 
\begin{align}\label{3.5}
I(\lambda) = 4if_0(i\ln(-i\lambda)), \qquad \lambda \in \C.
\end{align}
It is shown in \cite{FLhelmholtz} that $f_0(\varphi)$ admits an analytic continuation to complex values of $\varphi$ and hence the right-hand side of (\ref{3.5}) is well defined.

We recall that if the modified Helmholtz equation (\ref{2.1}) is defined in the {\it interior} domain $D_B$, then $W$ satisfies the equation
\begin{align}\label{3.6}
\int_{\partial D_B} W(z, \bar{z}, \lambda) = 0, \qquad \lambda \in \C.
\end{align}
It turns out that if $D_B$ is the interior of a polygon, then equation (\ref{3.6}) provides an efficient way for computing the generalized Dirichlet to Neumann map (see \cite{FF2011}). Indeed, for a well posed problem, one specifies a relationship between the Dirichlet and the Neumann boundary values. Thus, if the relevant polygon consists of $n$ sides, then equation (\ref{3.6}) contains $n$ unknowns. For example, in the case of the Dirichlet problem, equation (\ref{3.6}) contains the $n$ unknown Neumann boundary values. However, equation (\ref{3.6}) is valid for {\it all} complex values of $\lambda$. Thus, by expanding the unknown boundary values in an appropriate basis, and by evaluating equation (\ref{3.6}) at a sufficiently large number of collocation points, we can obtain the unknown coefficients appearing in the above expansions. Fornberg and colleagues have shown that by choosing Legendre polynomials as base functions, and by `overdetermining' the relevant linear system (i.e. by choosing the number of collocation points much larger than the number of unknown coefficients), the above technique yields a linear system with a small condition number. 
Furthermore, it is shown in \cite{FAS2013}, that by choosing the collocation points to be on certain curves, it is possible to obtain a condition number which is {\it independent} of $\beta$.

Equation (\ref{3.1}) differs from equation (\ref{3.6}) only by the presence of the explicit term defined in (\ref{3.5}). Thus, we can still apply the techniques developed for analyzing equation (\ref{3.6}), but we must supplement equation (\ref{3.1}) with an {\it additional} equation in order to compensate for the existence of the additional {\it unknown} function $f_0$.

\section{A Supplement to the Global Relation}\nequation
Suppose that $D$ is the exterior of the $n$-gon with vertices at $\{z_j\}_1^n$. 

It is shown in \cite{Sthesis} and \cite{FLhelmholtz} that in this case the following relation is valid:
\begin{align}\label{4.1}
\sum_{j=1}^n \int_{l_j} e^{i\beta(\lambda z - \frac{\bar{z}}{\lambda})} \hat{u}_j(\lambda) \frac{d\lambda}{\lambda} = 0, \qquad z \in D_B,
\end{align}
where $D_B$ denotes the interior of the polygon, the rays $l_j$ oriented from the origin to infinity are defined by 
\begin{align}\label{4.2}
l_j = \{\lambda \in \CÊ\, | \, \arg \lambda = -\arg(z_{j+1} - z_j)\}, \qquad j = 1, \dots, n,
\end{align}
and the functions $\{\hat{u}_j(\lambda)\}_1^n$ are defined by
\begin{align}\label{4.3}
\hat{u}_j(\lambda) = \int_{z_j}^{z_{j+1}} e^{-i\beta(\lambda z - \frac{\bar{z}}{\lambda})}\bigg[i u_{\mathcal{N}} + i\beta \bigg(\frac{1}{\lambda} \frac{d\bar{z}}{dt} + \lambda \frac{dz}{dt}\bigg)u\bigg] dt,
\end{align}
where $u_{\mathcal{N}}$ denotes the Neumann boundary value.
Let us parametrize the side $[z_j, z_{j+1}]$ by
\begin{align}\label{4.4}
z(t) = m_j + th_j, \qquad -1 \leq t \leq 1, 
\end{align}
where
\begin{align}\label{4.5}
m_j = \frac{z_j + z_{j+1}}{2}, \qquad h_j = \frac{z_{j+1} - z_j}{2}.
\end{align}
Then, the expression $\hat{u}_j(\lambda)$ defined in (\ref{4.3}) can be written in the form
\begin{align}\label{4.6}
\hat{u}_j(\lambda) = \hat{u}_j^{\mathcal{D}}(\lambda) + \hat{u}_j^{\mathcal{N}}(\lambda), \qquad j = 1, \dots, n,
\end{align}
where
\begin{align}\label{4.7}
  \hat{u}_j^{\mathcal{D}}(\lambda) & = i\beta\bigg(\frac{1}{\lambda} \bar{h}_j + \lambda h_j\bigg) e^{-i\beta(\lambda m_j - \frac{\bar{m}_j}{\lambda})} \int_{-1}^1 e^{-i\beta t(\lambda h_j - \frac{\bar{h}_j}{\lambda})} u dt,
  	\\ \label{4.8}
\hat{u}_j^{\mathcal{N}}(\lambda) & = i e^{-i\beta(\lambda m_j - \frac{\bar{m}_j}{\lambda})} \int_{-1}^1 e^{-i\beta t(\lambda h_j - \frac{\bar{h}_j}{\lambda})} u_{\mathcal{N}} dt.
\end{align}

\section{The Dirichlet Problem}\nequation
For the Dirichlet problem the functions $\{\hat{u}_j^{\mathcal{D}}(\lambda)\}_1^n$ are known (see (\ref{4.7})). In order to compute $\hat{u}_j^{\mathcal{N}}(\lambda)$, we approximate the Neumann boundary value on the side $j$ as follows:
\begin{align}\label{5.1}
  u_j^{\mathcal{N}}(t) \approx \sum_{m=0}^M c_m^{(j)} P_m(t), \qquad j = 1, \dots, n,
\end{align}
where $P_m$ is the Legendre polynomial of degree $m$, and the unknown constants $c_m^{(j)}$ are to be determined.

Equation (\ref{4.8}) yields
\begin{align}\label{5.2}
\hat{u}_j^{\mathcal{N}}(\lambda) = ie^{-i\beta(\lambda m_j - \frac{\bar{m}_j}{\lambda})} \sum_{m = 0}^M c_m^{(j)} \int_{-1}^1 e^{-i\beta t(\lambda h_j - \frac{\bar{h}_j}{\lambda})} P_m(t) dt. 
\end{align}
Thus, using the expression \cite{FAS2013}
\begin{align}\label{5.3}
  \int_{-1}^1 e^{-i\Lambda x} P_m(x) dx = \sum_{p = 0}^m \frac{A_{mp} e^{i\Lambda} + B_{mp}e^{-i\Lambda}}{\Lambda^{p+1}}, \qquad m \geq 0, \quad \Lambda \in \C \setminus \{0\},
\end{align}
where
\begin{align}\label{5.4}
A_{mp} = \frac{(m+p)!(-1)^{m+p}}{2^pp!(m-p)! i^{p+1}}, \qquad
B_{mp} = -\frac{(m+p)!}{2^pp!(m-p)!i^{p+1}},
\end{align}
equation (\ref{5.2}) becomes
\begin{align}\label{5.5}
\hat{u}_j^{\mathcal{N}}(\lambda) = ie^{-i\beta(\lambda m_j - \frac{\bar{m}_j}{\lambda})} \sum_{m=0}^M c_m^{(j)} \sum_{p=0}^m \frac{A_{mp} e^{i\Lambda} + B_{mp} e^{-i\Lambda}}{\Lambda^{p+1}}\bigg|_{\Lambda = \beta(\lambda h_j - \frac{\bar{h}_j}{\lambda})}.
\end{align}
Thus, equation (\ref{4.1}) yields
\begin{align}\label{5.6}
u^{\mathcal{D}}(z) + i \sum_{j = 1}^n \sum_{m=0}^M c_m^{(j)} \sum_{p=0}^m \frac{A_{mp} Q_{pj}^+(z) + B_{mp} Q_{pj}^-(z)}{\beta^{p+1}} = 0, \qquad z \in D_B,
\end{align}
where
\begin{align}\label{5.7}
u^{\mathcal{D}}(z) = \sum_{j=1}^n \int_{\hat{l}_j} e^{i\beta (\lambda z - \frac{\bar{z}}{\lambda})} \hat{u}_j^{\mathcal{D}}(\lambda) \frac{d\lambda}{\lambda},
\end{align}
and the functions $Q_{pj}^\pm(z)$ are defined by
\begin{align}\label{5.8}
Q_{pj}^\pm(z) = \int_{\hat{l}_j} \frac{e^{i\beta(\lambda z - \frac{\bar{z}}{\lambda} - \lambda m_j + \frac{\bar{m}_j}{\lambda} \pm (\lambda h_j - \frac{\bar{h}_j}{\lambda}))}}{(\lambda h_j - \frac{\bar{h}_j}{\lambda})^{p+1}} \frac{d\lambda}{\lambda},
\end{align}
with $\hat{l}_j$ denoting the contour $l_j$ with small indentations inserted appropriately so that $\hat{l}_j$ passes around the points $\lambda = \pm \exp(-i\arg h_j)$ (these points are {\it removable} singularities, thus we deform $l_j$ to $\hat{l}_j$ {\it before} splitting the relevant integrals).

Suppose we take the limit as $z$ approaches the $j_0$'th side. Then, using for this side the parametrization
\begin{align}\label{5.9}
  z_0 = m_{j_0} + t_0 h_{j_0}, \qquad -1 \leq t_0 \leq 1,
\end{align}
we find
\begin{align}\label{5.10}
& Q_{pj_0}^\sigma(z_0) = \int_{\hat{l}_{j_0}} 
 \frac{e^{i\beta (t_0 + \sigma)[\lambda h_{j_0} - \frac{\bar{h}_{j_0}}{\lambda}]}}{(\lambda h_{j_0} - \frac{\bar{h}_{j_0}}{\lambda})^{p+1}} \frac{d\lambda}{\lambda}, \qquad  \sigma = \pm 1,
 	\\\label{5.11}
& Q_{pj}^\sigma(z_0) = \int_{\hat{l}_j} 
 \frac{e^{i\beta[\lambda(m_{j_0} - m_j) + \lambda(t_0h_{j_0} + \sigma h_j) - \frac{\bar{m}_{j_0} - \bar{m}_j}{\lambda} - \frac{t_0\bar{h}_{j_0} + \sigma \bar{h}_j}{\lambda}]}}{(\lambda h_j - \frac{\bar{h}_j}{\lambda})^{p+1}} \frac{d\lambda}{\lambda}, \qquad j \neq j_0, \quad \sigma = \pm 1.
\end{align}
By employing the above limiting procedure, we can obtain $n$ additional equations. 
%However, we only need to supplement the global relation (\ref{3.1}) with {\it one} of these equations.

\section{The Exterior of a Square}
\begin{figure}
\begin{center}
 \begin{overpic}[width=.45\textwidth]{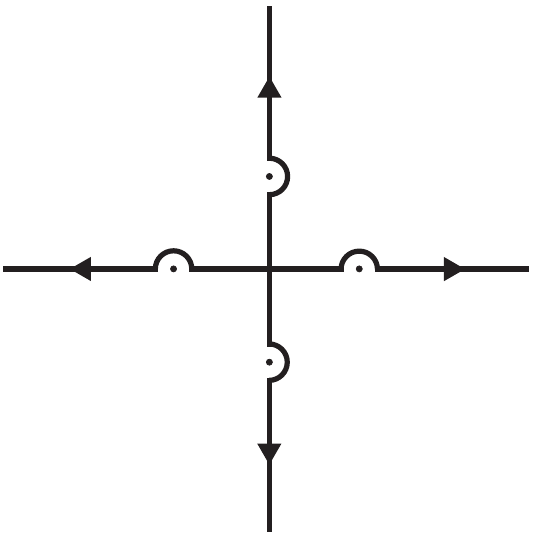}
  \put(65,43){$1$}
  \put(45,66){$i$}
  \put(27,43){$-1$}
  \put(40,31){$-i$}
  \put(82,40){$\hat{l}_1$}
  \put(42,81){$\hat{l}_2$}
  \put(14,40){$\hat{l}_3$}
  \put(42,14){$\hat{l}_4$}
    \end{overpic}
    \qquad \qquad
     \begin{figuretext}\label{Lhats.pdf}
       The indented contours $\hat{l}_j$, $j = 1, \dots, 4$, utilized for the exterior of the square. 
         \end{figuretext}
     \end{center}
\end{figure}

\nequation
Let $D_B$ be the square with corners at
$$z_1 = 1 +i, \qquad z_2 = -1 + i, \qquad z_3 = -1-i, \qquad z_4 = 1 - i.$$
Then, the four midpoints of the sides are
$$m_1 = i, \qquad m_2 = -1, \qquad m_3 = -i, \qquad m_4 = 1,$$
with the corresponding tangential vectors
$$h_1 = -1, \qquad h_2 = -i, \qquad h_3 = 1, \qquad h_4 = i.$$
Also,
$$l_1 = (0, -\infty), \qquad l_2 = (0, i\infty), \qquad l_3 = (0, \infty), \qquad l_4 = (0, -i\infty).$$
Thus, if $j_0 = 1$ so that $z_0 = i - t_0$, $-1 \leq t_0 \leq 1$, belongs to the top side, then
\begin{align*}
& Q_{p1}^\sigma(z_0) = \int_{\hat{l}_1} 
 \frac{e^{-i\beta (t_0 + \sigma)[\lambda - \frac{1}{\lambda}]}}{(-\lambda + \frac{1}{\lambda})^{p+1}} \frac{d\lambda}{\lambda}, && \sigma = \pm 1,
 	\\
& Q_{p2}^\sigma(z_0) = \int_{\hat{l}_2}  \frac{e^{i\beta[i(1-\sigma)(\lambda + \frac{1}{\lambda}) + (1-t_0)(\lambda - \frac{1}{\lambda})]}}{(-i\lambda - \frac{i}{\lambda})^{p+1}} \frac{d\lambda}{\lambda}, && \sigma = \pm 1,
 	\\
& Q_{p3}^\sigma(z_0) = \int_{\hat{l}_3}  \frac{e^{i\beta[2i(\lambda + \frac{1}{\lambda}) + (\sigma - t_0)(\lambda - \frac{1}{\lambda})]}}{(\lambda - \frac{1}{\lambda})^{p+1}} \frac{d\lambda}{\lambda}, && \sigma = \pm 1,
 	\\
& Q_{p4}^\sigma(z_0) = \int_{\hat{l}_4}  \frac{e^{i\beta[i(1+\sigma)(\lambda + \frac{1}{\lambda}) -(1+t_0)(\lambda - \frac{1}{\lambda})]}}{(\lambda i + \frac{i}{\lambda})^{p+1}} \frac{d\lambda}{\lambda}, && \sigma = \pm 1.
\end{align*}
We choose the indented contours $\hat{l}_j$, $j = 1, \dots, 4$, displayed in Figure \ref{Lhats.pdf}.

The exponential $e^{-i\beta (t_0 + \sigma)[\lambda - \frac{1}{\lambda}]}$
has decay for $\im \lambda < 0$ if $\sigma = 1$ and for $\im \lambda > 0$ if $\sigma = -1$. 
Using that
$$\underset{\lambda = -1}{\res}\frac{e^{-i\beta (t_0 + \sigma)[\lambda - \frac{1}{\lambda}]}}{(-\lambda + \frac{1}{\lambda})^{p+1}} \frac{1}{\lambda}
= \frac{1}{p!} \frac{d^p}{d\lambda^p}\bigg|_{\lambda = -1} \bigg(\frac{\lambda}{1-\lambda}\bigg)^{p+1} e^{-i\beta (t_0 + \sigma)[\lambda - \frac{1}{\lambda}]} \frac{1}{\lambda},$$
we find
\begin{align*}
 Q_{p1}^\sigma(z_0) = &\; \int_0^{e^{-\frac{3i\pi}{4}}\infty} 
 \frac{e^{-i\beta (t_0 + \sigma)[\lambda - \frac{1}{\lambda}]}}{(-\lambda + \frac{1}{\lambda})^{p+1}} \frac{d\lambda}{\lambda}
 	\\
& + \frac{2\pi i}{p!} \frac{d^p}{d\lambda^p}\bigg|_{\lambda = -1} \bigg(\frac{\lambda}{1 - \lambda}\bigg)^{p+1} e^{-i\beta (t_0 + \sigma)[\lambda - \frac{1}{\lambda}]} \frac{1}{\lambda}, \qquad \sigma = 1, 
\end{align*}
and
\begin{align*}
Q_{p1}^\sigma(z_0) = \int_0^{e^{\frac{3i\pi}{4}}\infty}
 \frac{e^{-i\beta (t_0 + \sigma)[\lambda - \frac{1}{\lambda}]}}{(-\lambda + \frac{1}{\lambda})^{p+1}} \frac{d\lambda}{\lambda}, \qquad \sigma = -1.
\end{align*}
Similarly, we find
\begin{align*}
Q_{p2}^\sigma(z_0) = & \; \int_0^{\infty e^{\frac{i\pi}{4}}}  \frac{e^{i\beta[i(1-\sigma)(\lambda + \frac{1}{\lambda}) + (1-t_0)(\lambda - \frac{1}{\lambda})]}}{(-i\lambda - \frac{i}{\lambda})^{p+1}} \frac{d\lambda}{\lambda}, \qquad \sigma = \pm 1,
	\\
 Q_{p3}^\sigma(z_0) = &\; \int_0^{\infty e^{\frac{i\pi}{4}}}  \frac{e^{i\beta[2i(\lambda + \frac{1}{\lambda}) + (\sigma - t_0)(\lambda - \frac{1}{\lambda})]}}{(\lambda - \frac{1}{\lambda})^{p+1}} \frac{d\lambda}{\lambda}, \qquad \sigma = 1, 
	\\
Q_{p3}^\sigma(z_0) = &\; \int_0^{\infty e^{-\frac{i\pi}{4}}}  \frac{e^{i\beta[2i(\lambda + \frac{1}{\lambda}) + (\sigma - t_0)(\lambda - \frac{1}{\lambda})]}}{(\lambda - \frac{1}{\lambda})^{p+1}} \frac{d\lambda}{\lambda}
	\\
&- \frac{2\pi i }{p!} \frac{d^p}{d\lambda^p}\bigg|_{\lambda = 1} 
\bigg(\frac{\lambda}{\lambda+1}\bigg)^{p+1} e^{i\beta[2i(\lambda + \frac{1}{\lambda}) + (\sigma - t_0)(\lambda - \frac{1}{\lambda})]} \frac{1}{\lambda}, \qquad \sigma = -1,
	\\
Q_{p4}^\sigma(z_0) = &\; \int_0^{\infty e^{-\frac{i\pi}{4}}}  \frac{e^{i\beta[-(1+t_0)(\lambda - \frac{1}{\lambda}) + i(1+\sigma)(\lambda + \frac{1}{\lambda})]}}{(\lambda i + \frac{i}{\lambda})^{p+1}} \frac{d\lambda}{\lambda}, \qquad \sigma = \pm 1.
\end{align*}
The above integrals have exponential decay and are suitable for numerical evaluation.
Similar considerations lead to analogous expressions for $Q_{pj}^\pm(z_0)$ for $z_0$ belonging to the left, bottom, and right sides of the square. 

Substituting the above expressions into relation (\ref{5.6}) evaluated at $z = z_0$ where $z_0$ belongs to one of the sides of the square, we find relations for the $c_m^{(j)}$'s which can be used to supplement the exterior global relation. Taking into consideration that the global relation now involves {\it one} additional unknown (the scattering amplitude $f_0$), one would expect that it is sufficient to supplement the global relation with the equation obtained by varying $z_0$ on any {\it one} of the sides of the square. However, it appears that in order to obtain accurate results it is necessary to consider all four sides. 

It was noted earlier that the term $u^{\mathcal{D}}$ defined by (\ref{5.7}) is known. However, the question of computing this term numerically in the most efficient way remains open. 
In the examples below we compute the known expression (\ref{5.7}) by approximating the given Dirichlet data in terms of Legendre polynomials.

\subsection{Examples}
Let $D$ be the domain exterior to the square with vertices at $\pm 1 \pm i$. The function
\begin{align}\label{uHn}
u(r, \varphi) = H_n^{(1)}(2i\beta r) (A e^{in\varphi} +  B e^{-in\varphi}),
\end{align}
satisfies the modified Helmholtz equation (\ref{2.1}) in $D$ for any $n \geq 1$ and any complex constants $A$ and $B$. 

\begin{figure}
\begin{center}
\bigskip\bigskip
 \begin{overpic}[width=.42\textwidth]{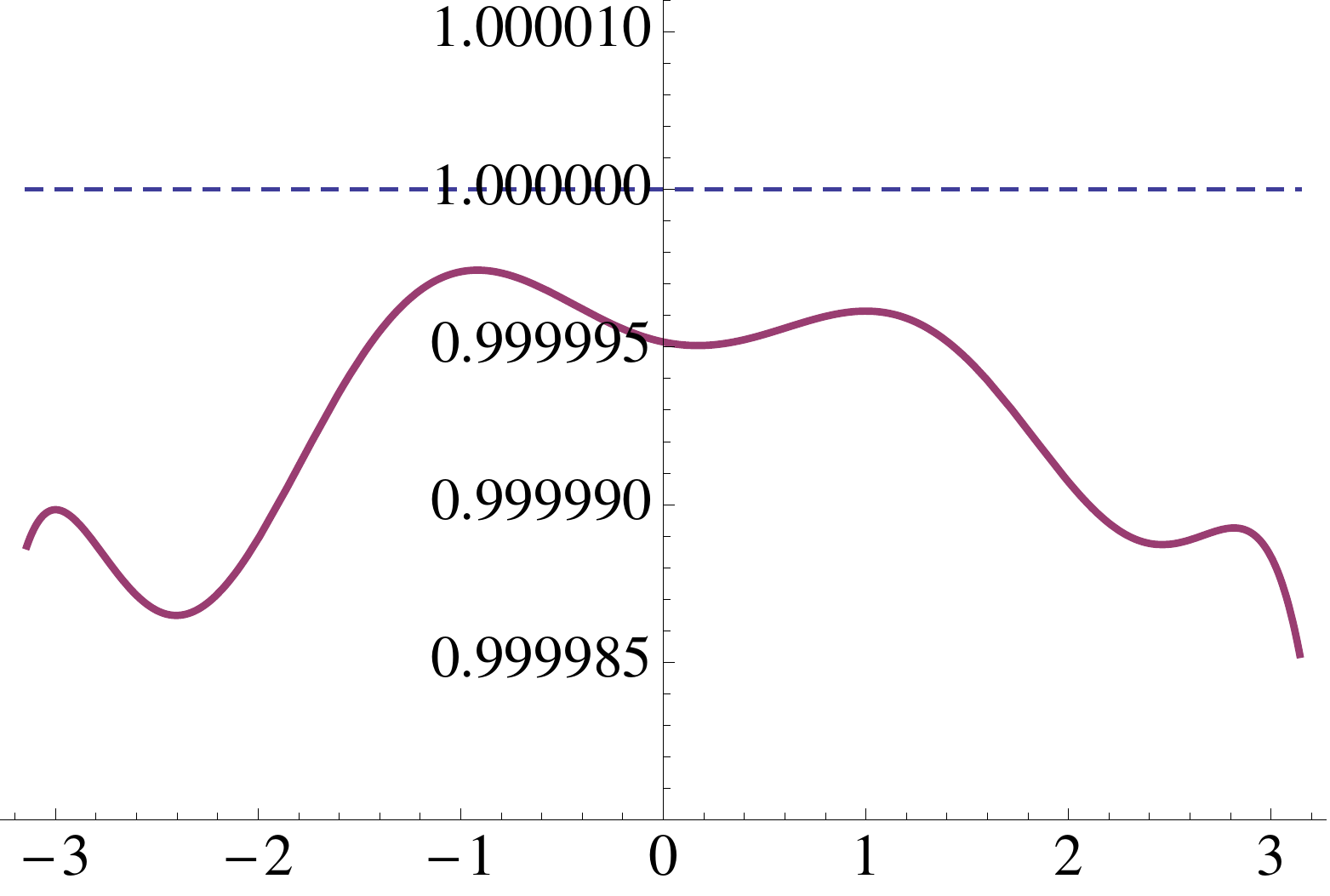}
  \put(37,74){$\re f_0(\varphi)$}
  \put(103,4){$\varphi$}
    \end{overpic}
    \qquad\quad
 \begin{overpic}[width=.42\textwidth]{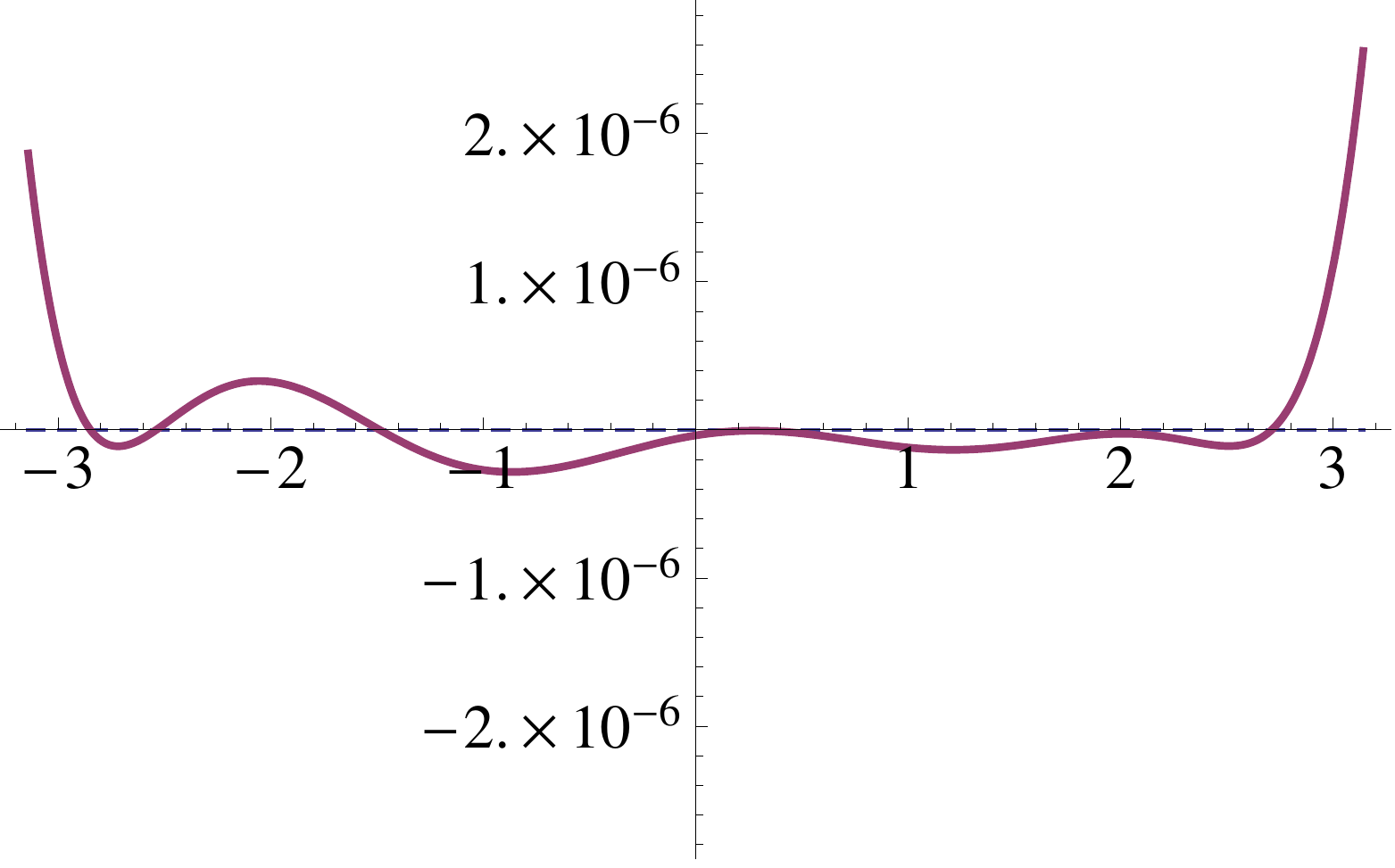}
  \put(37,68){$\im f_0(\varphi)$}
   \put(103,30){$\varphi$}
   \end{overpic}
   \bigskip
     \begin{figuretext}\label{H0fig1}
       The real and imaginary parts of the scattering amplitude $f_0(\varphi)$ for Example 1. Numerical values (solid) compared with true values (dashed).
         \end{figuretext}
     \end{center}
\end{figure}
\begin{figure}
\begin{center}
\bigskip\bigskip
 \begin{overpic}[width=.42\textwidth]{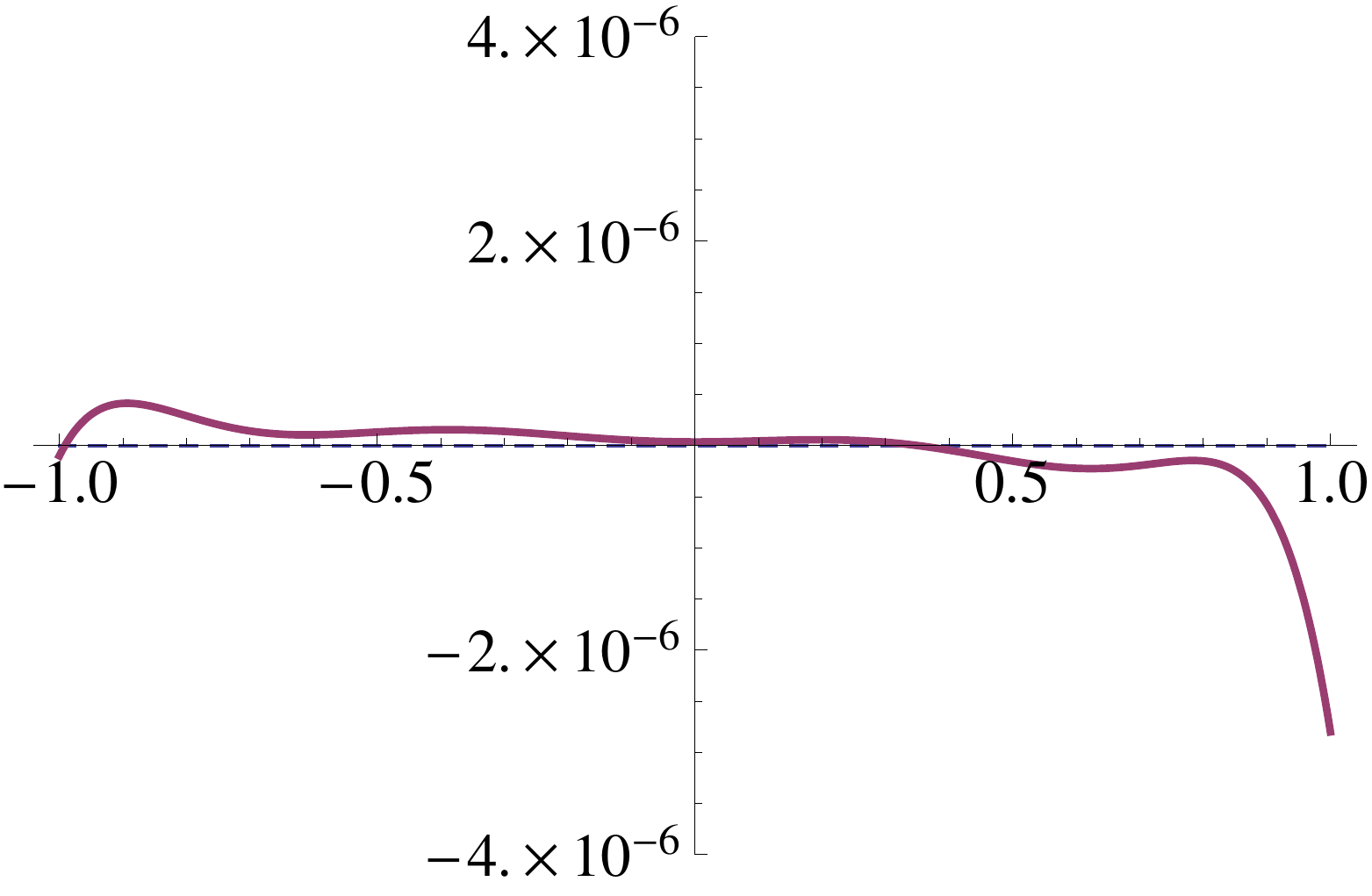}
  \put(30,72){$\re u_y(i-t)$}
  \put(103,31){$t$}
    \end{overpic}
    \qquad\quad
 \begin{overpic}[width=.42\textwidth]{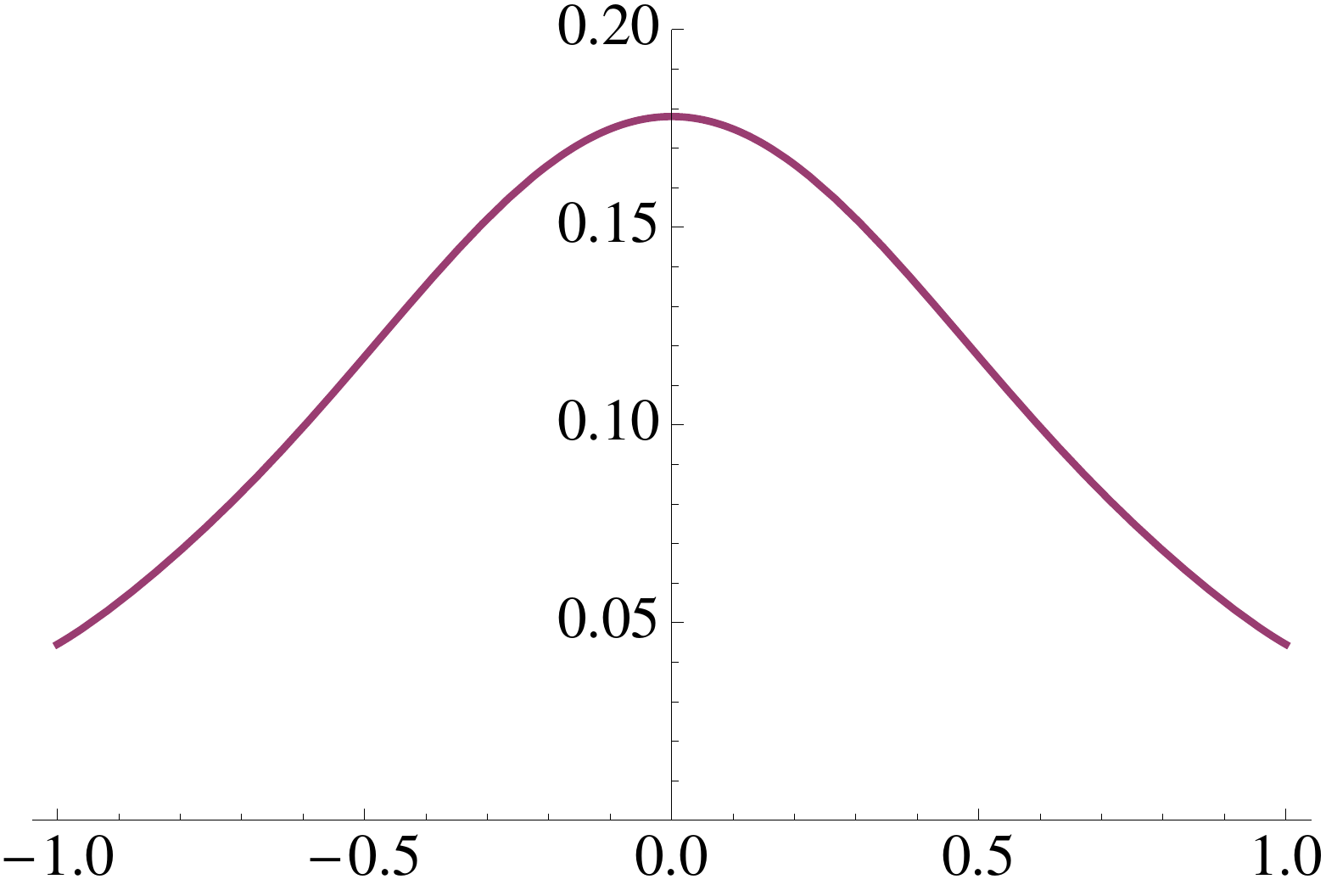}
  \put(30,72){$\im u_y(i-t)$}
  \put(103,4){$t$}
    \end{overpic}
\bigskip
     \begin{figuretext}\label{H0fig2}
       The real and imaginary parts of the Neumann value on the top side for Example 1. 
       %The error is too small to be visible for the imaginary part. 
         \end{figuretext}
     \end{center}
\end{figure}

\subsubsection{Example 1}
For the first example, we let $u(r, \varphi) = H_0^{(1)}(2i \beta  r)$. This solution satisfies
$$u(r,\varphi) = \sqrt{\frac{1}{\pi i \beta  r}} e^{i(2i\beta r - \frac{\pi}{4})} \biggl(1-\frac{1}{16 \beta  r} + \frac{9}{512 \beta ^2 r^2} + O(r^{-3})\biggr), \qquad r \to \infty.$$
Hence equation (\ref{3.4}) is satisfied with $f_0(\varphi) \equiv 1$. 
%Hence, by (\ref{3.5}), $I(\lambda) = 4i$.
We use the values of $u$ on the four sides of the square as the given Dirichlet data. Letting $\beta = 1$, we apply the above technique with $M = 8$ and with $54$ collocation points. The numerical scheme yields the scattering amplitude $f_0(\varphi)$ displayed in Figure \ref{H0fig1}. The numerically obtained and true Neumann values on the top side are shown in Figure \ref{H0fig2}; the results for the other three sides are similar. The numerics gives the correct profiles with a maximal error of approximately $10^{-4}$.

\subsubsection{Example 2}
\begin{figure}
\begin{center}
\bigskip\bigskip
 \begin{overpic}[width=.42\textwidth]{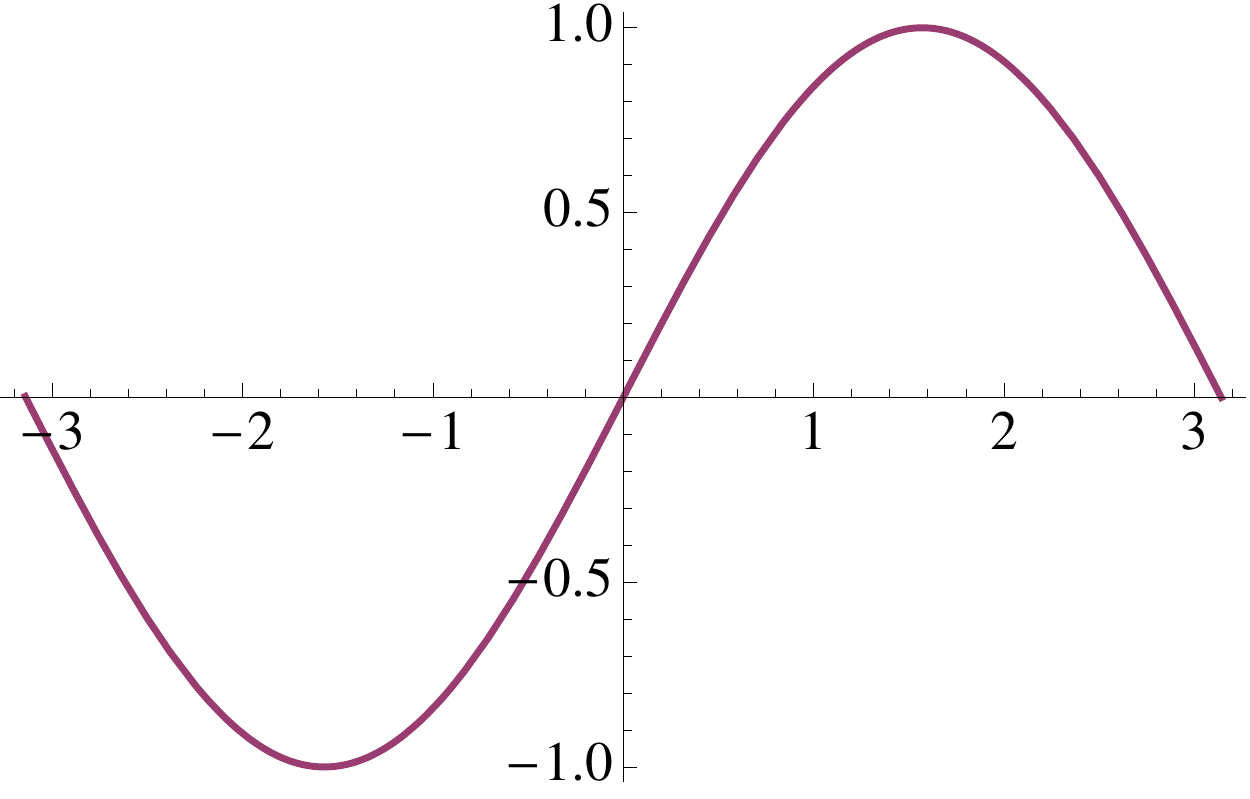}
  \put(37,70){$\re f_0(\varphi)$}
  \put(103,31){$\varphi$}
    \end{overpic}
    \qquad\quad
 \begin{overpic}[width=.42\textwidth]{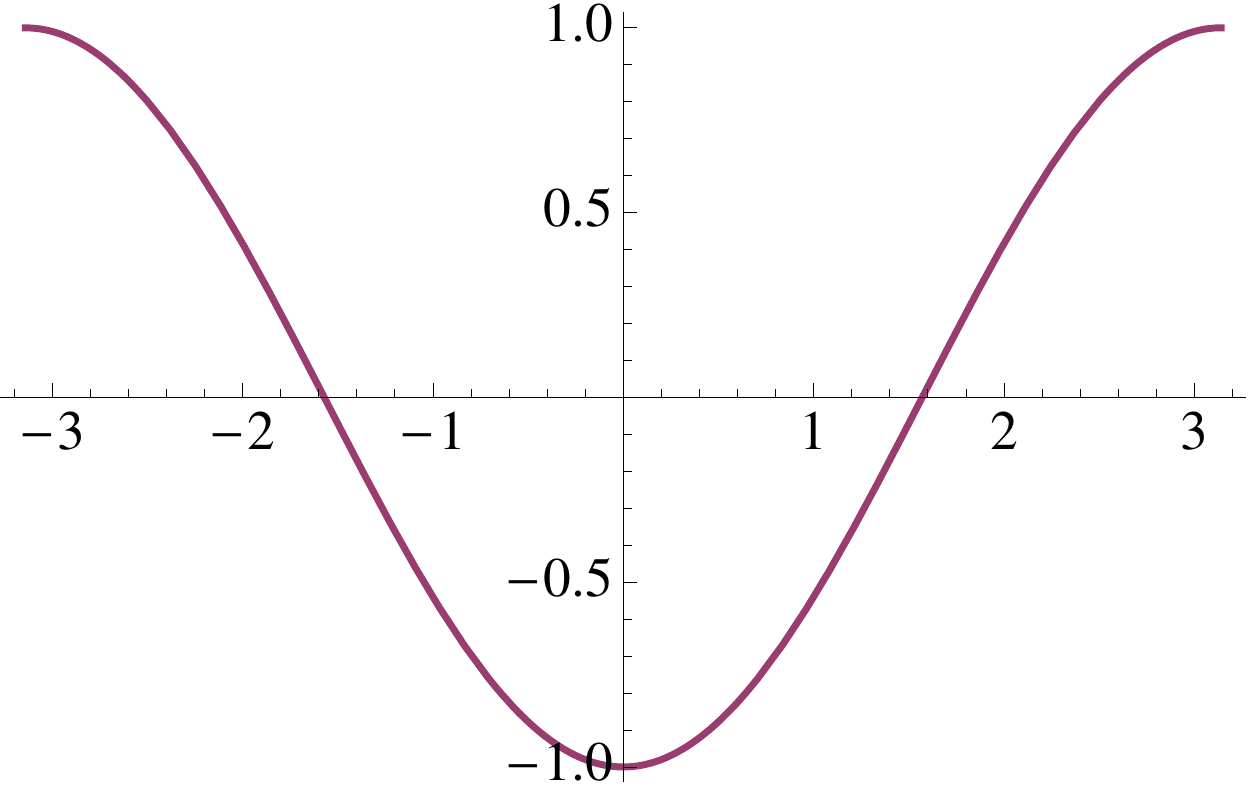}
  \put(37,70){$\im f_0(\varphi)$}
   \put(103,31){$\varphi$}
   \end{overpic}
   \bigskip
     \begin{figuretext}\label{H1fig1}
       The real and imaginary parts of the scattering amplitude $f_0(\varphi)$ for Example 2. 
       \end{figuretext}
     \end{center}
\end{figure}

\begin{figure}
\begin{center}
\bigskip\bigskip
 \begin{overpic}[width=.42\textwidth]{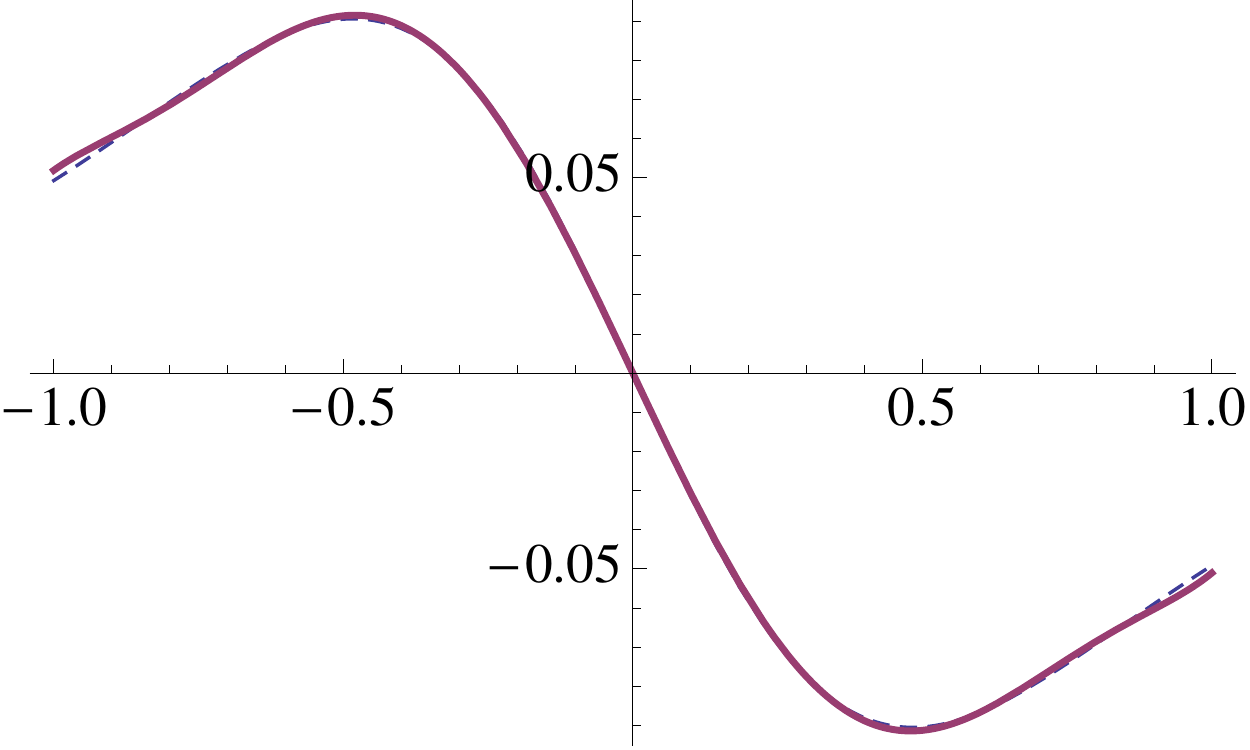}
  \put(32,68){$\re u_y(i-t)$}
  \put(103,28){$t$}
    \end{overpic}
    \qquad\quad
 \begin{overpic}[width=.42\textwidth]{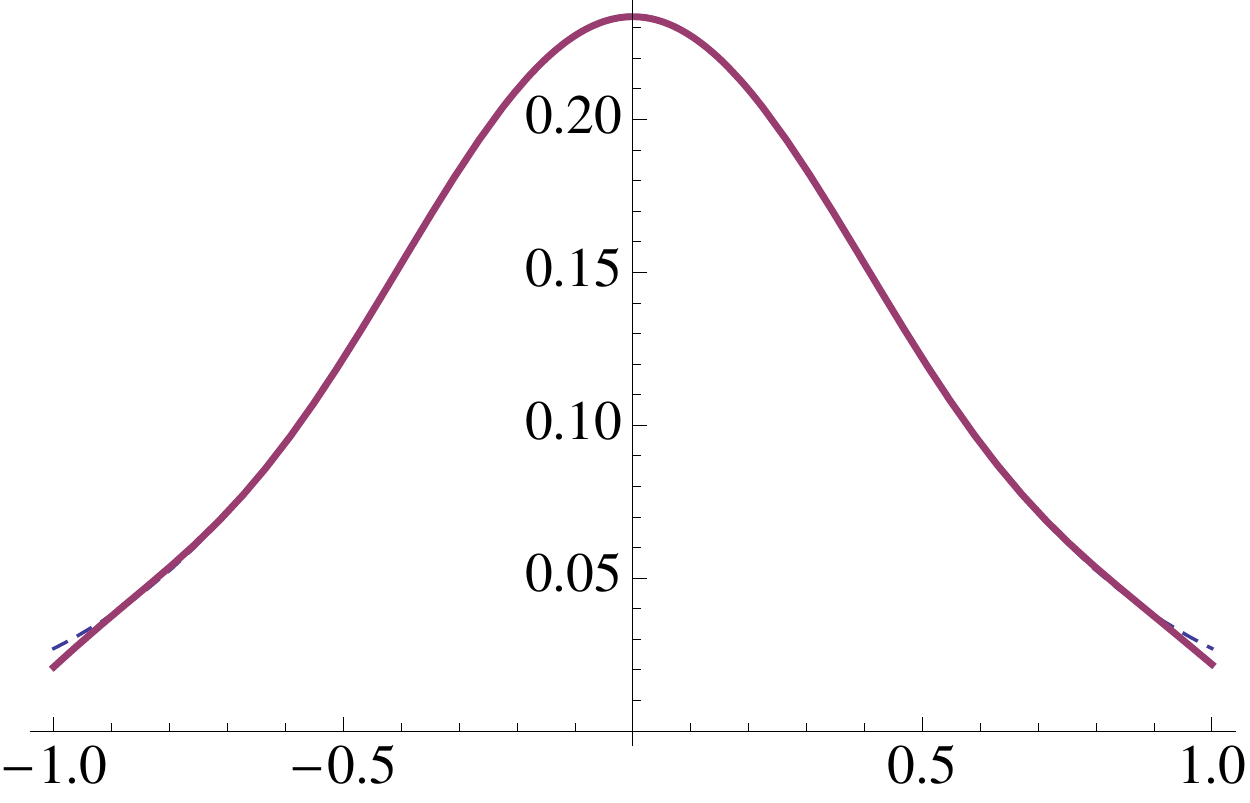}
  \put(32,70){$\im u_y(i-t)$}
  \put(103,4){$t$}
    \end{overpic}
\bigskip
     \begin{figuretext}\label{H1fig2}
       The real and imaginary parts of the Neumann value on the top side for Example 2. 
         \end{figuretext}
     \end{center}
\end{figure}

For the second example, we let $u(r, \varphi) = H_1^{(1)}(2i \beta  r)e^{i\varphi}$. This solution satisfies
$$u(r,\varphi) = \sqrt{\frac{1}{\pi i\beta  r}} e^{i(2i\beta r - \frac{\pi}{4})} \biggl(-ie^{i\varphi} +\frac{3 e^{i\varphi }}{16 i\beta  r} + \frac{15 i e^{i \varphi}}{512 \beta ^2 r^2} + O(r^{-3})\biggr), \qquad r \to \infty.$$
Hence equation (\ref{3.4}) is satisfied with $f_0(\varphi) = -ie^{i\varphi}$. 
%Hence, by (\ref{3.5}), $I(\lambda) = 4i/\lambda$.
Letting $\beta = 1$, the numerical scheme with $M = 8$ and with $54$ collocation points yields the scattering amplitude $f_0(\varphi)$ displayed in Figure \ref{H1fig1}. The numerically obtained and true Neumann values on the top side are shown in Figure \ref{H1fig2}. The numerics gives the correct profiles with a maximal error of approximately $0.002$.


\section{Conclusions}\nequation
We have illustrated a novel numerical technique for linear elliptic PDEs formulated in the {\it exterior} of a polygon. This technique, which is based on the Unified Transform, uses two main ingredients: (a) the global relation; in the case of the modified Helmholtz equation, the global relation is given by (\ref{3.1}), where $I(\lambda)$ is defined by equation (\ref{3.5}). This equation can be analyzed numerically in the same way that the global relation is analyzed in the case of linear elliptic PDEs formulated in the {\it interior} of a polygon. However, now the global relation, in addition to the unknown boundary values, also involves the unknown function $f_0$ (the scattering amplitude). Hence, it must be supplemented with an additional equation. (b) Additional equations can be obtained by taking the limit of a certain equation which is valid for $z$ in the interior of the given polygon, as $z$ approaches the boundary; in the case of the modified Helmholtz equation, this equation is given by (\ref{5.6}). Actually, by taking the limit as $z$ approaches each side of an $n$-gon, one obtains $n$ additional equations. These equations provide the formulation in the complex Fourier plane (complex extension of $\lambda$) of the classical boundary integral method. By exploiting the analyticity structure in the complex Fourier plane, it is possible to perform contour deformations and to obtain integrals which decay exponentially for large $\lambda$. This yields a most efficient numerical evaluation of these integrals. Detailed comparisons will be presented elsewhere.

\bigskip
\noindent
{\bf Acknowledgement} {\it The authors acknowledge support from the EPSRC, UK.}

\bibliographystyle{plain}
\bibliography{is}

\end{document}